\documentclass[a4paper,10pt]{amsart}
\usepackage{amssymb}
\usepackage{cite}	
\usepackage{tikz}
\usepackage[left=2cm,right=2cm,top=3cm,bottom=2cm]{geometry} % page settings
\usepackage{amsmath} % provides many mathematical environments & tools
\allowdisplaybreaks
\usepackage{amsfonts}
\usepackage{amsmath}
\usepackage{amssymb}
\usepackage{systeme}
\usepackage{mathrsfs}
\usepackage[colorlinks=true]{hyperref}
\numberwithin{equation}{section}
\newtheorem{theorem}{Theorem}[section]
\newtheorem{prop}[theorem]{Proposition}
\newtheorem{lm}[theorem]{Lemma}

\newtheorem{defi}[theorem]{Definition}

\newtheorem{rem}{Remarks}
\newcommand{\R}{\mathbb{R}}

\newcommand{\N}{\mathbb{N}^\ast}

\newcommand{\supp}{\operatorname{supp}}

\usepackage{setspace}

\def\ds{\displaystyle}

\def\dive{\mathrm{ div}}
\newenvironment{preuve}{{\noindent{\textbf{\large Proof.}}}}{\hfill {\rule{2.5mm}{2.5mm}}}
\makeatletter
\let\@msm@th@eqref\eqref
\renewcommand{\eqref}[1]{%
	\begingroup
	\leavevmode
	\color{red}%
	\hypersetup{linkbordercolor=[named]{red}}%
	\@msm@th@eqref{#1}%
	\endgroup
}
\makeatother
\makeatletter
\newcommand{\leqnomode}{\tagsleft@true\let\veqno\@@leqno}
\newcommand{\reqnomode}{\tagsleft@false\let\veqno\@@eqno}
\makeatother
\usepackage[stable]{footmisc}
\author[M.~Amara]{Mustapha Amara${}^\ast$\protect\footnotetext{${}^\ast$Corresponding author}}
\author[C.~Katar]{Chaala Katar}
\author[M.~Ltifi]{Maroua Ltifi}
\address{Department of Mathematics, Faculty of Science of Gab\`es, Research Laboratory Mathematics and Applications LR17ES11; Tunisia}
\email{\sl Mostafa.Amara@fsg.u-gabes.tn}
\email{\sl katarchaala123@gmail.com}
\email{\sl widaltifi@gmail.com}

\title[Generalization of $3D-NSE$ Global Weak  Solution with damping ]
{Generalization of $3D-NSE$ Global Weak  Solution with damping}

\begin{document}
	\begin{abstract}
		In  this paper, we prove the global existence, uniqueness and the continuity in $L^2$ of the incompressible Navier-Stokes equations with damping $f(|u|)u$, where $f$ is an increasing, convex and differentiable function on $\R^+$, null at zero. 
	\end{abstract}
	%@@@@@@@@@@@@@@@@@@@@@@@@@@@@@@@@@@@@%@@@@@@@@@@@@@@@@@@@@@@@@@@@@@@@@@@@@%@@@@@@@@@@@@@@
	
	\subjclass[2020]{35-XX, 35Q30, 76D05, 76N10}
	\keywords{Navier-Stokes equations, Friedrich method, global weak solution}

	%@@@@@@@@@@@@@@@@@@@@@@@@@@@@@@@@@@@@%@@@@@@@@@@@@@@@@@@@@@@@@@@@@@@@@@@@@%@@@@@@@@@@@@@@
	\maketitle
	%@@@@@@@@@@@@@@@@@@@@@@@@@@@@@@@@@@@@%@@@@@@@@@@@@@@@@@@@@@@@@@@@@@@@@@@@@%@@@@@@@@@@@@@@@
	\section{\bf Introduction}
	In this paper, we consider the Cauchy problem of three-dimensional $(3D)$ incompressible Navier-Stokes equations with the additional nonlinear term $f(|u|)u$	introduced in the following equation
	\leqnomode
	\begin{equation*}\tag*{$(NSD)$}\label{NSD}
		\begin{cases}
			\partial_t u
			-\nu\Delta u+ u\cdot\nabla u  + f(|u|)u =\;\;-\nabla \pi&\mbox{ in } \R^+\times \R^3,\\
			\dive(u) = 0& \mbox{ in } \R^+\times \R^3,\\
			u(0,x) =u^0(x) &\mbox{ in }\R^3,
		\end{cases}
	\end{equation*}
	\reqnomode
	where $u=u(t,x)=(u_1(t,x),u_2(t,x),u_3(t,x))$ and $\pi=\pi(t,x)$ represent the unknown velocity vector and the unknown pressure scalar of the fluid at the point $(t,x)\in \R^+\times \R^3$, respectively, and $\nu>0$ denotes the fluid's viscosity. The term $f(|u|)u$ serves as the damping term, with $f$ being a positive function defined on $\R^+$. This term models various physical phenomena, including flow through porous media, drag or friction effects, and certain dissipative mechanisms (see \cite{HP,HL} and reference therein).\\
	The nonlinear term  $\left(u\cdot\nabla u\right):=u_1\partial_1u+u_2\partial_2u+u_3\partial_3u$, the fact that $\dive(u)=0$, can be written in the following form:
	\begin{equation*}
		\left(u\cdot\nabla u\right)=\dive(u\otimes u)=\left(\dive(u_1u),\dive(u_2u),\dive(u_3u)\right).
	\end{equation*}
	While $u^0=u^0(x)=(u_1^0(x),u_2^0(x),u_3^0(x))$ is an initial given velocity.\\
	
	The global existence of weak solutions for the initial value problem of the classical incompressible Navier-Stokes equations ( the case where $f$ is identically zero) was established long ago by Leray and Hopf (see \cite{L,H}). However, the uniqueness of these weak solutions remains an open question in dimensions $d\geq 3$.\\
	
	The damping term, in the mathematical framework, is added to control the nonlinear term to solve the uniqueness problem. The equation \ref{NSD} is a generalization of the incompressible Navier-Stokes equations with polynomial damping, when $f(x)=\alpha x^{\beta-1}$, $\alpha>0$  and $\beta\geq 1$, which is studied in \cite{CJ} by Cia and Jui, where they proved the global existence of the weak solution in 
	\begin{equation}
		L^{\infty}(\R^{+},L^{2}(\R^{3}))\cap L^{2}(\R^{+},\dot{H}^{1}(\R^{3}))\cap L^{\beta+1}(\R^{+}\times\R^{3}).
	\end{equation}
	The uniqueness and continuity of this weak solution have recently been proved in \cite{BB} by Blel and Benameur when $\beta>3$ and $\alpha>0$.\\
	
	In 2022, Benameur in \cite{JBB} introduced a new damping term when $f(x)=\alpha(e^{\beta x^2}-1)$, and established the global existence of weak solutions when $\alpha,\beta>0$, in
	\begin{equation}
		L^{\infty}(\R^{+},L^{2}(\R^{3}))\cap L^{2}(\R^{+},\dot{H}^{1}(\R^{3}))\cap L^{p}(\R^{+}\times\R^{3})\cap\mathcal{E}_\beta,
	\end{equation}
	for all $p\geq 4$, where $\mathcal{E}_\beta=\left\{u:\R^+\times \R^3\rightarrow \R\ \mbox{ measurable and }\ (e^{\beta|u|^2}-1)|u|^2\in L^1(\R^+\times\R^3)\right\}.$ \\
	For the uniqueness and continuity of the weak solution, Blel and Benameur in \cite{BB1} indicate that for the damping term when $f(x)=\alpha(e^{\beta x^r}-1)$, $\alpha,\beta>0$ and $r\geq 1$, by adapting the same proof of result in \cite{JBB}, we have the existence of the global weak solution in 
	\begin{equation}
		C_b(\R^{+},L^{2}(\R^{3}))\cap L^{2}(\R^{+},\dot{H}^{1}(\R^{3})) \cap\mathcal{E}_\beta^r,
	\end{equation}
	where $\mathcal{E}_\beta^r=\left\{u:\R^+\times \R^3\rightarrow \R\ \mbox{ measurable and }\ (e^{\beta|u|^r}-1)|u|^2\in L^1(\R^+\times\R^3)\right\}.$ Moreover, in the same paper, they show that for the term damping when $f(x)=\alpha(e^{\beta x^4}-1)$, $\alpha,\beta>0$, we have the uniqueness and the continuity in $L^2(\R^3).$\\
	
	The purpose of this paper is to study the well-posedness of the incompressible Navier-Stokes equations with a damping term, where $f$ is an increasing, convex and differentiable function on $\R^+$, null at zero. We will prove the existence of global weak solutions for the Cauchy problem \ref{NSD} which verified the following estimates:
	\begin{equation*}
		\|u(t)\|^{2}_{L^{2}}+2\int_{0}^{t}\|\nabla u(z)\|^{2}_{L^{2}}dz+2\int_{0}^{t}\|f(|u(z)|)|u(z)|^2\|_{L^1}dz\leq\|u^{0}\|^2_{L^{2}},
	\end{equation*} 
	to guarantee that solution belong to
	$$L^{\infty}(\R^{+},L^{2}(\R^{3}))\cap L^{2}(\R^{+},\dot{H}^{1}(\R^{3}))\cap \mathcal{H},$$
	where 
	$$\mathcal H:=\{u: \R^+\times\R^3\rightarrow\R^3\ \mbox{ measurable and }f(|u|)|u|^{2}\in L^{1}(\R^{+}\times \R^3)\}.$$ 
	Furthermore, we will show that if there exists $c>0$ and $p>2$ such that 
	\begin{equation}\label{Condition}
		f(x)\geq cx^p,\quad \forall x\geq 0,
	\end{equation} 
	then the weak solution is unique and belong to $$C_b(\R^+,L^2(\R^3))\cap L^{p+2}(\R^+\times  \R^3).$$
	
	Before dealing with global existence, we give the definition of the weak solution of \ref{NSD} system
	\begin{defi}
		The couple of functions $(u,\pi)$ is called a weak solution of problem \ref{NSD} if, for all $T > 0$, the following conditions are satisfied:
		\begin{enumerate}
			\item[\textbf{1.}] $u\in L^\infty([0,T],L^2(\R^3))\cap L^2([0,T],\dot{H}^1(\R^3))$ and $f(|u|)|u|^2\in L^1([0,T]\times \R^3).$
			\item[\textbf{2.}]  The couple $(u,\pi)$  satisfies 
			\begin{equation}
				\partial_{t} u-\nu\Delta u+u\cdot \nabla u+f(|u|)u=-\nabla \pi,\quad \mbox{ in }D'([0,T]\times \R^3).
			\end{equation}
			In other words, for any $\phi\in D([0,T]\times \R^3)$ such that $\dive(\phi)(t,x)=0$, for all $(t,x)\in [0,T]\times \R^3$, and $\phi(T)=0$, we have 
			\begin{equation}
				\int_{0}^{T}(u,\partial_{t} \phi)_{L^2}=-(u^0,\phi(0))_{L^2}+	\nu\int_{0}^{T}(\nabla u,\nabla \phi)_{L^2}+	\int_{0}^{T}((u\cdot\nabla)u, \phi)_{L^2}+	\int_{0}^{T}(f(|u|)u, \phi)_{L^2}.
			\end{equation}
			\item[\textbf{3.}]  $\dive(u)(x,t)=0$ for a.e $(t,x)\in[0,T]\times\R^{3}.$
		\end{enumerate}
	\end{defi}
	\begin{theorem}\label{th1} Let $f$ be an increasing convex differentiable function on $\R^+$ such that $f(0)=0$. If $u^0\in L^{2}(\R^{3})$ such that $\dive(u^0)=0$, then, the system \ref{NSD} have a global weak solution $(u,\pi)$ such that
		\begin{equation}
			u\in L^{\infty}(\R^{+},L^{2}(\R^{3}))\cap L^{2}(\R^{+},\dot{H}^{1}(\R^{3}))\cap C(\R^+,H^{-2}(\R^3))\cap\mathcal H,
		\end{equation}
		satisfied, for all $t>0,$
		\begin{equation}\label{1.8}
			\|u(t)\|^{2}_{L^{2}}+2\int_{0}^{t}\|\nabla u(z)\|^{2}_{L^{2}}dz+2\int_{0}^{t}\int_{\R^3}f(|u(z,x)|)|u(z,x)|^2dxdz\leq\|u^{0}\|^2_{L^{2}}.
		\end{equation} 
		Moreover, if there exist $c>0$ and $p>2$ such that for all $x\in \R^+$
		\begin{equation}\label{1.4}
			f(x)\geq cx^p,\quad \forall x\geq 0,
		\end{equation}
		then the solution $u$ is unique and 
		\begin{equation}
			u\in C_b(\R^+,L^2(\R^3))\cap L^{p+2}(\R^+\times \R^3).
		\end{equation}
	\end{theorem}
	\begin{rem}
		\begin{enumerate}
			\item[] 
			\item[\textbf{1.}] If $(u,\pi)$ is a solution from Theorem \ref{th1}, then $f(|u|)|u|\in L^1_{loc}(\R^+,L^1(\R^3)),$ moreover, formally,
			$$\pi=(-\Delta)^{-1}\left[\dive(u\cdot \nabla u)+\dive(f(|u|)u)\right]\in L^1_{\mbox{loc}}(\R^+,H^{-s}(\R^3)),\quad \forall s>\frac{3}{2}.$$
			Indeed: For $T>0$, consider the subset $\mathcal{A}_T=\left\{(t,x)\in [0,T]\times \R^3;\ |u(t,x)|\leq 1\right\}$, then
			\begin{align*}
				\int_{[0,T]\times \R^3}f(|u|)|u|&=\int_{\mathcal{A}_T}f(|u|)|u|+\int_{\mathcal{A}_T^c}f(|u|)|u|.
			\end{align*}
			Almost for everything $(t,x)\in\mathcal{A}_T$ we have $|u(t,x)|\in [0,1]$ and the fact that $f$ is differentiable and convex function, according to the Lemma \ref{Lemma2.4}, we get
			\begin{equation*}
				f(|u(t,x)|)\leq f'(1)|u(t,x)|.
			\end{equation*}
			Almost for everything $(t,x)\in\mathcal{A}_T^c$ we have $|u(t,x)|\geq 1\Rightarrow  |u(t,x)|\leq |u(t,x)|^2.$\\
			Therefore
			\begin{align*}
				\int_{[0,T]\times \R^3}f(|u|)|u|&\leq \int_{\mathcal{A}_T}f'(1)|u|^2+\int_{\mathcal{A}_T^c}f(|u|)|u|^2\\
				&\leq f'(1)T\|u\|_{L^\infty{(\R^+,L^2)}}^2+\|f(|u|)|u|^2\|_{L^1(\R^+\times \R^3)}<\infty.
			\end{align*}
			In the other hand, if $s> \frac{3}{2}$, we have 
			\begin{align*}
				\|(-\Delta)^{-1}(\dive(u\cdot \nabla u)(t)\|_{H^{-s}}&\leq\|(u\otimes u)(t)\|_{H^{-s}}\\
				&\leq C_s \|\mathcal{F}(u\otimes u)(t)\|_{L^\infty}&C_s=\left(\int_{\R^3}\dfrac{d\xi}{(1+|\xi|^2)^{s}}\right)^{\frac{1}{2}}\\
				&\leq C_s\|(u\otimes u)(t)\|_{L^1}\\
				&\leq C_s \|u(t)\|^{2}_{L^2}\in L_{\mbox{loc}}^1(\R^+),
			\end{align*}
			and 
			\begin{align*}
				\|(-\Delta)^{-1}(\dive(f(|u|) u)(t)\|_{H^{-s}}&\leq\left(\int_{\R^3}\frac{1}{|\xi|^2(1+|\xi|^2)^s}|\mathcal{F}(f(|u|) u)(t,\xi)|^2d\xi\right)^{\frac{1}{2}}\\
				&\leq \tilde{C_s}\|\mathcal{F}(f(|u|) u)(t)\|_{L^\infty} &\hskip-2cm\tilde{C_s}=\left(\int_{\R^3}\dfrac{d\xi}{|\xi|^2(1+|\xi|^2)^{s}}\right)^{\frac{1}{2}}\\
				&\leq \tilde{C_s}\|f(|u(t)|) u(t)\|_{L^1}\in L_{\mbox{loc}}^1(\R^+).
			\end{align*}
			Which implies that $\pi\in L^1_{\mbox{loc}}(\R^+,H^{-s}(\R^3)),\quad \forall s>\frac{3}{2}.$
			\item[\textbf{2.}] The fact that $u\in L^\infty(\R^+,L^2(\R^3))\cap C(\R^+,H^{-2}(\R^3))$ implies that $u\in C_b(\R^+,H^{-r}(\R^3))$, for all $r>0$.
			\item [\textbf{3.}] If $ f $ satisfies \eqref{1.4}, then, for all $ x \geq 0 $,  
			$$
			f(x) \geq cx^{p} \implies f(x)x^2 \geq cx^{p+2}.
			$$
			Thus,  
			$$
			\int_{\R^+ \times \R^3} |u|^{p+2} \leq \frac{1}{c} \int_{\R^+ \times \R^3} f(|u|)|u|^2 < \infty,
			$$
			which implies that $ u \in L^{p+2}(\R^+ \times \R^3) $.
			\item[\textbf{4.}] If $f(x)=\alpha x^{\beta-1}$, $\alpha>0$ and $\beta\geq 2$, then $f$ is increasing convex differentiable function on $\R^+$ and $f(0)=0$.
			Moreover, if $\beta>3$ then $f$ verifies the condition \eqref{1.4} $(p=\beta-1\mbox{ and } c=\alpha)$, which means that there exists a unique weak global solution 
			$$u\in C_b(\R^+,L^2(\R^3))\cap L^{\beta+1}(\R^+\times \R^3).$$
			\item[\textbf{5.}] If $f(x)=\alpha (e^{\beta x^r}-1)$, $\alpha,\beta>0$ and $r\geq 1$, then $f$ is increasing convex differentiable function on $\R^+$ and $f(0)=0$.
			Moreover $f$ verifies the condition \eqref{1.4} $\left(p=r+1\mbox{ and } c=\frac{\alpha\beta^{\frac{r}{r+1}}}{r+1}\right)$, which means that there exists a unique weak global solution 
			$$u\in C_b(\R^+,L^2(\R^3))\cap L^{r+3}(\R^+\times \R^3).$$
		\end{enumerate}
	\end{rem}	
	
	The remainder of this paper is organized as follows. In Section 2, we introduce the notations, definitions, and preliminary results. Section 3 is dedicated to the proof of our main result: the global existence, uniqueness, and continuity of the solution to the Cauchy problem \ref{NSD} in $ L^2(\R^3) $.
	
	\section{\bf Preliminary result}
	Let us present the elementary information that is used in this paper:
	\begin{enumerate}
		\item [\textbf{1.}] $\mathcal{S}'(\R^3)$ is the space of tempered distributions.
		\item[\textbf{2.}] The Fourier transform and its inverse are defined by
		$$
		\mathcal{F}(f)(\xi)=\widehat{f}(\xi):=\int_{\R^3}e^{-ix\cdot\xi}f(x)dx,
		$$
		and
		$$
		\mathcal{F}^{-1}(g)(x):=(2\pi)^{-3}\int_{\R^3}e^{i\xi\cdot x}g(\xi)d\xi,
		$$
		respectively.
		\item[\textbf{3.}] The convolution product of a suitable pair of function $f$ and $g$ on $\R^3$ is given by
		$$
		(f\ast g)(x):=\int_{\R^3}f(y)g(x-y)dy.
		$$
		\item[\textbf{4.}] The tensor product is given by
		$$
		f\otimes g:=(g_1f,g_2f,g_3f),
		$$
		where $f=(f_1,f_2,f_3)$ and $g=(g_1,g_2,g_3)\in S'(\R^3).$ Moreover
		$$
		\dive(f\otimes g):=(\dive(g_1f),\dive(g_2f),\dive(g_3f)),
		$$
		and if $\dive(g)=0$ we get
		$$
		\dive(f\otimes g):=g_1\partial_1f+g_2\partial_2f+g_3\partial_3f:=(g\cdot\nabla) f.
		$$
		\item[\textbf{5.}] Let $(X,\|\cdot\|_X)$ be a normed space. We define, for all $T>0$
		$$ L^\infty([0,T],X)=\left\{u:[0,T]\rightarrow X\mbox{ measurable and }\|\|u(t)\|_X\|_{L^\infty([0,T])}<\infty\right\},$$
		and 
		$$ L^\infty(\R^+,X)=\left\{u:\R^+\rightarrow X\mbox{ measurable and }\|\|u(t)\|_X\|_{L^\infty(\R^+)}<\infty\right\}.$$
		The $L^\infty([0,T],X)-$norm (resp. $ L^\infty(\R^+,X)-$ norm) is given by
		$$\|u\|_{L^\infty([0,T],X)}:=\|\|u(t)\|_X\|_{L^\infty([0,T])},\quad \left(\mbox{resp. } \|u\|_{L^\infty(\R^+,X)}:=\|\|u(t)\|_X\|_{L^\infty(\R^+)}\right).$$
		\item[\textbf{6.}] Let $(X,\|\cdot\|_X)$ be a normed space. We define, for all $T>0$ and $p\geq 1$
		$$ L^p([0,T],X)=\left\{u:[0,T]\rightarrow X\mbox{ measurable and }\int_{(0,T)}\|u(t)\|_X^pdt<\infty\right\},$$
		and 
		$$ L^p(\R^+,X)=\left\{u:\R^+\rightarrow X\mbox{ measurable and }\int_{\R^+}\|u(t)\|_X^p<\infty\right\}.$$
		The $L^p([0,T],X)-$norm (resp. $ L^p(\R^+,X)-$ norm ) is given by
		$$\|u\|_{L^p([0,T],X)}:=\left(\int_{(0,T)}\|u(t)\|_X^pdt\right)^{\frac{1}{p}},\quad \left(\mbox{resp.} \|u\|_{L^p(\R^+,X)}:=\left(\int_{\R^+}\|u(t)\|_X^pdt\right)^{\frac{1}{p}}\right).$$
		\item[\textbf{7.}] Let $s\in \R$. Define the non-homogeneous Sobolev space by 
		$$H^s(\R^3)=\{u\in \mathcal S'(\R^3):\ \widehat{u}\in L^1_{loc}(\R^3) \mbox{ and }\int_{\R^3}(1+|\xi|^2)^{s}|\widehat{u}(\xi)|^2d\xi<\infty\},$$
		where $H^s(\R^3)$ is endowed with the norm
		$$\|u\|_{H^s}:=\left(\int_{\R^3}(1+|\xi|^2)^{s}|\widehat{u}(\xi)|^2d\xi\right)^{\frac{1}{2}},\quad\forall u\in H^s(\R^3),$$
		and, moreover, $H^s(\R^3)$ is endowed with the inner product
		$$(u,v)_{H^s}:=\int_{\R^3}(1+|\xi|^2)^{s}\widehat{u}(\xi)\overline{\widehat{v}(\xi)}d\xi,\quad\forall u,v\in H^s(\R^3).$$
		\item[\textbf{8.}] Let $s\in \R$. Define the homogeneous Sobolev space by 
		$$\dot{H}^s(\R^3)=\{u\in \mathcal S'(\R^3):\ \widehat{u}\in L^1_{loc}(\R^3) \mbox{ and }\int_{\R^3}|\xi|^{2s}|\widehat{u}(\xi)|^2d\xi<\infty\},$$
		where $\dot{H}^s(\R^3)$ is endowed with the norm
		$$\|u\|_{\dot{H}^s}:=\left(\int_{\R^3}|\xi|^{2s}|\widehat{u}(\xi)|^2d\xi\right)^{\frac{1}{2}},\quad\forall u\in H^s(\R^3),$$
		and, moreover, $\dot{H}^s(\R^3)$ is endowed with the inner product
		$$(u,v)_{\dot{H}^s}:=\int_{\R^3}|\xi|^{2s}\widehat{u}(\xi)\overline{\widehat{v}(\xi)}d\xi,\quad\forall u,v\in \dot{H}^s(\R^3).$$
		\item[\textbf{9.}] Let $R>0$. Define the Friedrich operator $J_R$ by
		$$J_R(u)=\mathcal F^{-1}\left(\xi\mapsto\mathcal{X}_{B(0,R)}(\xi)\mathcal{F}({u})(\xi)\right),$$
		where $B(0,R)$ is the ball in $\R^3$ of centre $0$ and radius $R$.
		\item[\textbf{10.}] The Leray projector $\mathbb P:(L^2(\R^3))^3\rightarrow L^2_\sigma(\R^3)$ is defined by
		$$\mathcal F(\mathbb P (u))(\xi)=\widehat{u}(\xi)-(\widehat{u}(\xi),\xi)_{\R^3}\frac{\xi}{|\xi|^2}=M(\xi)\widehat{u}(\xi),$$
		where $L^2_\sigma(\R^3)=\{f\in (L^2(\R^3))^3:\; \dive(f)=0\}$ and $M(\xi)$ is the matrix $\left(\delta_{i,j}-\dfrac{\xi_i\xi_j}{|\xi|^2}\right)_{1\leq i,j\leq 3}$.\\ Define also the
		operator $A_R(D)$ on $(L^2(\R^3))^3 $ by
		$$A_R(D)u:=J_R(\mathbb{P} (u))=\mathcal{F}^{-1}\left(\xi \mapsto M(\xi)\mathcal{X}_{B(0,R)}(\xi)\mathcal{F}({u})(\xi)\right).$$	
	\end{enumerate}
	Now, we need to establish some results that are key points in the development of the techniques applied in this work. The reader is referred to \cite{HBAF,HB,JYC} for more details of the proofs.
	\begin{prop}\label{prop1} Let $H$ be a Hilbert space.
		\begin{enumerate}
			\item If $(x_n)_n$ is a bounded sequence of elements in $H$, then there is a sub-sequence $(x_{{\varphi(n)}})_n$ such that
			$$(x_{{\varphi(n)}}|y)\rightarrow (x|y),\quad\forall y\in H.$$
			\item If $x\in H$ and $(x_n)_n$ is a bounded sequence of elements in $H$ such that
			$$(x_n|y)\rightarrow (x|y),\quad\forall y\in H.$$
			Then $\|x\|\leq\displaystyle\liminf_{n\rightarrow\infty}\|x_n\|.$
			\item If $x\in H$ and $(x_n)_n$ is a bounded sequence of elements in $H$ such that $\displaystyle\limsup_{n\rightarrow\infty}\|x_n\|\leq \|x\|$ and
			$$
			(x_n|y)\rightarrow (x|y),\quad\forall y\in H,$$
			then $\displaystyle\lim_{n\rightarrow\infty}\|x_n-x\|=0.$
		\end{enumerate}
	\end{prop}
	\begin{lm}\label{LP}
		Let $s_1,\ s_2$ be two real numbers such that $s_1+s_2>0$. If $s_1<3/2$, then, there exists a constant  $C=C(s_1,s_2,d)>0$, such that: for all $u,v\in \dot{H}^{s_1}(\R^3)\cap \dot{H}^{s_2}(\R^3)$, we have $uv \in \dot{H}^{s_1+s_2-\frac{3}{2}}(\R^3)$ and
		$$\|uv\|_{\dot{H}^{s_1+s_2-\frac{3}{2}}}\leq C (\|u\|_{\dot{H}^{s_1}}\|v\|_{\dot{H}^{s_2}}+\|u\|_{\dot{H}^{s_2}}\|v\|_{\dot{H}^{s_1}}).$$
		Moreover, if also $s_2<3/2$, then, for all $u \in \dot{H}^{s_1}(\R^3)$ and $v\in\dot{H}^{s_2}(\R^3)$, then  $uv \in \dot{H}^{s_1+s_2-\frac{3}{2}}(\R^3)$ and
		$$\|fg\|_{\dot{H}^{s_1+s_2-\frac{3}{2}}}\leq C \|f\|_{\dot{H}^{s_1}}\|g\|_{\dot{H}^{s_2}}.$$
	\end{lm}
	\begin{lm}\label{lm45}
		Let $a,b\in \R$ such that $a<b$ and let $g,h:[a,b]\rightarrow\R^+$ two continuous functions and $C>0$ such that
		\begin{equation*}\label{LG2}\forall t\in[a,b],\quad g(t)\leq C+\int_a^th(z)g(z)dz.
		\end{equation*}
		Then $$\forall t\in[a,b],\quad g(t)\leq C\exp\left[\int_a^th(z)dz\right].$$
	\end{lm}
	%%%%%%%%%%%%%%%%%%%
	\begin{lm}\label{Lemma2.4}
		Let $ R > 0 $, and let $ f $ be a convex and differentiable function on $[0, R]$, such that $f(0)=0$. Then, for any $ x, y \in [0, R] $,
		$$|f(x)-f(y)|\leq f'(R)|x-y|.$$
		In particular, for any $x\in [0,R]$,
		$$f(x)\leq f'(R)x.$$
	\end{lm}
	\begin{preuve}
		All we need to do is to apply the finite growth theorem to $f$ on $[0,R]$ and use the fact that $f'$ is increasing on $[0,R]$.
	\end{preuve}
	\begin{lm}\label{Lemma2.5}
		Let $f:\R^+\rightarrow \R^+$ a function and let $c>0$ and $p>0$ such that
		\begin{equation}
			f(x)\geq cx^p,\quad x\geq 0,
		\end{equation}
		then, for all $u,v\in \R^3$, we have
		\begin{equation}
			\left(f(|u|)u-f(|v|)v,u-v\right)_{\R^3}\geq \frac{c}{4}\left(|u|^{p}+|v|^{p}\right)|u-v|^2.
		\end{equation}
	\end{lm}
	\begin{preuve}
		Let $u,v\in \R^3$. Then,
		\begin{align*}
			\left(f(|u|)u-f(|v|)v,u-v\right)_{\R^3}
			&=  f(|u|)|u|^2+ f(|v|)|v|^2-(f(|u|)+f(|v|))\left(u,v\right)_{\R^3}\\
			&=  f(|u|)|u|^2+ f(|v|)|v|^2-\frac{1}{4}(f(|u|)+f(|v|))\left(|u+v|^2-|u-v|^2\right)\\
			&=\frac{1}{4}(f(|u|)+f(|v|))|u-v|^2+f(|u|)|u|^2+ f(|v|)|v|^2-\frac{1}{4}(f(|u|)+f(|v|))|u+v|^2.
		\end{align*}
		We now need to prove that
		\begin{align*}
			I=f(|u|)|u|^2+ f(|v|)|v|^2-\frac{1}{4}(f(|u|)+f(|v|))|u+v|^2\geq 0,
		\end{align*}
		Indeed, we have:
		\begin{align*}
			\frac{1}{4}(f(|u|)+f(|v|))|u+v|^2&\leq \frac{1}{2}(f(|u|)+f(|v|))\left(|u|^2+|v|^2\right)\\
			&\leq \frac{1}{2}f(|u|)|u|^2+\frac{1}{2}f(|u|)|v|^2+ \frac{1}{2}f(|v|)|u|^2+ \frac{1}{2}f(|v|)|v|^2.
		\end{align*}
		Thus,
		\begin{align*}
			I&\geq f(|u|)|u|^2+ f(|v|)|v|^2-\frac{1}{2}f(|u|)|u|^2-\frac{1}{2}f(|u|)|v|^2- \frac{1}{2}f(|v|)|u|^2- \frac{1}{2}f(|v|)|v|^2\\
			&\geq \frac{1}{2}f(|u|)|u|^2+\frac{1}{2} f(|v|)|v|^2-\frac{1}{2}f(|u|)|v|^2- \frac{1}{2}f(|v|)|u|^2\\
			&\geq \frac{1}{2}f(|u|)(|u|^2-|v|^2)-\frac{1}{2}f(|v|)(|u|^2-|v|^2)\\
			&\geq \frac{1}{2}\left(f(|u|)-f(|v|)\right)(|u|^2-|v|^2)\\
			&\geq 0.
		\end{align*}
		where we use the fact that $f$ is increasing function.\\
		
		Therefore,
		\begin{align*}
			\left(f(|u|)u-f(|v|)v,u-v\right)_{\R^3}&\geq \frac{1}{4}(f(|u|)+f(|v|))|u-v|^2
		\end{align*}
		Using the fact that
		$$f(x)\geq c x^{p},\quad \forall x\in \R^+,$$
		we get
		\begin{align*}
			\left(f(|u|)u-f(|v|)v,u-v\right)_{\R^3}&\geq \frac{c}{4}(|u|^p+|v|^p)|u-v|^2.
		\end{align*}
		Hence the result.
	\end{preuve}
	\section{\bf Proof of Theorem \ref{th1}}
	\subsection{Existence}~\\
	\textbf{Step 1: Approximate the system}.\\
	We consider the following approximate system
	\leqnomode
	\begin{equation}\tag*{$(NSD_n)$}
		\begin{cases}
			\partial_t u -\Delta J_{n}(u)+J_{n}\left(J_n(u)\cdot\nabla J_n(u)\right) +J_{n}\left(f(|J_{n}(u)|)J_{n}(u)\right)=-\nabla \pi_{n}&\mbox{ in }\R^+\times\R^3, \\
			\pi_n=J_n(-\Delta)^{-1}\Big[\dive\left(J_n(u)\cdot\nabla J_n(u)\right)+\dive\left(f(|J_{n}(u)|)J_{n}(u)\right)\Big]& \\
			\dive(u)=0&\mbox{ in }\R^+\times\R^3,\\
			u(0,x)=J_n(u^0(x))&\mbox{  in }\R^3.
		\end{cases}
	\end{equation}
	Cauchy-Lipschitz theorem gives a unique solution $u_n$ in the space $\mathcal{C}^1(\R^+,L^2(\R^{3}))$. Moreover $J_{n}u_n=u_{n}$, then 
	\begin{equation}
		\begin{cases}
			\partial_t u_n -\Delta u_{n}+A_n(D)\left(u_n\cdot\nabla u_n\right) +A_{n}(D)\left(f(|u_{n}|)u_{n}\right)=0&\mbox{ in }\R^+\times\R^3, \\
			\dive( u_n)=0&\mbox{ in }\R^+\times\R^3,\\
			u_n(0,x)=J_n(u^0(x))&\mbox{  in }\R^3.
		\end{cases}
	\end{equation} 
	\reqnomode
	Moreover, for all $t>0$
	\begin{equation}\label{Eq0}
		\|u_{n}(t)\|^{2}_{L^{2}}+2\int_{0}^{t}\|\nabla u_{n}(s)\|^{2}_{L^{2}}ds+2\int_{0}^{t}\|f(|u_n(s)|)|u_n(s)|^2\|_{L^1}ds \leq \|u^{0}\|^{2}_{L^{2}}.
	\end{equation}
	\textbf{Step 2: Passage to the limit}.\\
	$\bullet$ The sequence $(u_n)_n$ is equicontinuity in $C_b(\R^+,H^{-2}(\R^3))$.\\
	Indeed, for all $t_1,t_2\in \R^+$ such that $t_1\leq t_2$ we have 
	\begin{align*}
		\|u_n(t_2)-u_n(t_1)\|_{H^{-2}}&\leq \int_{t_1}^{t_2}\|\Delta u_n\|_{H^{-2}}+\int_{t_1}^{t_2}\|A_n(D)\left(u_n\cdot\nabla u_n\right)\|_{H^{-2}}+\int_{t_1}^{t_2}\|A_n(D)\left(f(|u_n|) u_n\right)\|_{H^{-2}}\\
		&\leq \int_{t_1}^{t_2}\| u_n\|_{L^2}+\int_{t_1}^{t_2}\|u_n\cdot\nabla u_n\|_{H^{-2}}+\int_{t_1}^{t_2}\|f(|u_n|) u_n\|_{H^{-2}}.
	\end{align*}
	We need to control the terms to the right of the inequality term by term, starting with the first term, we have 
	\begin{equation}
		\int_{t_1}^{t_2}\| u_n\|_{L^2}\leq \int_{t_1}^{t_2}\| u^0\|_{L^2}\leq \| u^0\|_{L^2}(t_2-t_1).
	\end{equation}
	For the second term, using Lemma \ref{LP} for $s_1=0$ and $s_2=1$ ($s_1+s_2-\frac{3}{2}=-\frac{1}{2}$) and using the fact that $\dive(u_n)=0$, we obtain
	$$\|u_n\cdot\nabla u_n\|_{H^{-2}}=\|\dive(u_n\otimes u_n)\|_{H^{-2}}\leq \|u_n\otimes u_n\|_{H^{-1}}\leq \|u_n\otimes u_n\|_{\dot{H}^{-\frac{1}{2}}}\leq C\|u_n\|_{L^2}\|\nabla u_n\|_{L^2}.$$
	Thus 
	\begin{align*}
		\int_{t_1}^{t_2}\|u_n\cdot\nabla u_n\|_{H^{-2}}&\leq C \int_{t_1}^{t_2}\| u_n\|_{L^2}\|\nabla u_n\|_{L^2}\\
		&\leq C \|u^0\|_{L^2}\int_{t_1}^{t_2}\|\nabla u_n\|_{L^2}\\
		&\leq C \|u^0\|_{L^2}\sqrt{t_2-t_1}\left(\int_{t_1}^{t_2}\|\nabla u_n\|^2_{L^2 }\right)^\frac{1}{2}\\
		&\leq C \|u^0\|^2_{L^2}\sqrt{t_2-t_1}.
	\end{align*}
	For the last term, for all $R>0$ we have
	\begin{align*}
		\|f(|u_n|) u_n\|_{H^{-2}}\leq C \|f(|u_n|) u_n\|_{L^1}&=C\int_{\R^3} f(|u_n|) |u_n| dx\\
		&=C\int_{A_{n,R,t}} f(|u_n|) |u_n| dx+C\int_{A_{n,R,t}^c} f(|u_n|) |u_n| dx.
	\end{align*}
	where 
	$$A_{n,R,t}:=\{x\in \R^3;\ |u_n(t,x)|\leq R\}.$$
	Using Lemma \ref{Lemma2.4}, we get for all $x\in A_{n,R,t}$
	$$f(|u_n|)\leq f'(R)|u_n|\Rightarrow f(|u_n|)|u_n|\leq f'(R)|u_n|^2.$$
	Moreover, for all $x\in  A_{n,R,t}^c$
	$$f(|u_n|)|u_n|\leq \frac{1}{R}f(|u_n|)|u_n|^2.$$
	Thus
	\begin{align*}
		\|f(|u_n|) u_n\|_{H^{-2}}\leq C \|f(|u_n|) u_n\|_{L^1}	&\leq C f'(R)\|u_n\|_{L^2}^2 dx+\dfrac{C}{R}\| f(|u_n|) |u_n|^2\|_{L^1}.
	\end{align*}
	As a result
	\begin{align*}
		\int_{t_1}^{t_2}\|f(|u_n|) u_n\|_{H^{-2}}&\leq C f'(R)	\int_{t_1}^{t_2}\|u_n\|^2 dx+\dfrac{C}{R}	\int_{t_1}^{t_2}\| f(|u_n|) |u_n|^2\|_{L^1}.\\
		&\leq  Cf'(R)\|u^0\|^2_{L^2}(t_2-t_1)+\dfrac{C}{R}\|u^0\|_{L^2}^2.
	\end{align*}
	Let $\varepsilon>0$ and take $$R>\dfrac{4C\|u^0\|^2_{L^2}}{\varepsilon}.$$
	We assume $$0<\delta\leq \dfrac{1}{4}\min\left\{\dfrac{\varepsilon}{\|u^0\|_{L^2}+1},\dfrac{\varepsilon^2}{C^2(\|u^0\|_{L^2}+1)^2},\dfrac{\varepsilon}{Cf'(R)(\|u^0\|_{L^2}+1)}\right\},$$ 
	then for any $t_1,t_2\in \R^+$ such that $|t_2-t_1|\leq \delta$ and for any $n\in \N$ we have
	$$\|u_n(t_2)-u_n(t_1)\|_{H^{-2}}\leq \| u^0\|_{L^2}|t_2-t_1|+C \|u^0\|^2_{L^2}\sqrt{|t_2-t_1|}+Cf'(R)\|u^0\|_{L^2}^2|t_2-t_1|+\dfrac{C}{R}\|u^0\|_{L^2}^2\leq \varepsilon.$$
	Which implies that $\{u_{n},\ n\in\N\}$ is equicontinuous  in $C_b(\R^+,H^{-2}(\R^3))$.\\
	
	\noindent$\bullet$ Let $ (T_q)_{q \in \N} \subset \R^{+,\ast}$ be a sequence such that $ T_q < T_{q+1} $ and $ T_q \to \infty $ as $ q \to \infty $, and for any $ q \in \N$, take $ \theta_q \in D(\R^3) $ a test function defined by:
	$$\begin{cases}
		\theta_q(x)=1,&\forall x\in B(0, q+\frac{5}{4})\\
		\theta_q(x)=0,& \forall x\in B(0, q+2)^c\\
		0\leq \theta_q\leq 1.	    
	\end{cases}
	$$
	Using the fact that $ \{ u_{n},\ n \in \N\} $ is equicontinuous in $C_b(\R^+,H^{-2}(\R^3))$, and by applying a classical argument combining Ascoli's theorem with the Cantor diagonal process, we obtain a non-decreasing sequence $ \varphi: \N\to \N $ and a function $ u \in L^\infty(\R^+, L^2(\R^3)) \cap C(\R^+, H^{-2}(\R^3)) $, such that for all $ q \in \N $, we have:
	\begin{equation}\label{Lim0}
		\lim_{n\rightarrow\infty}\|\theta_q(u_{{\varphi(n)}}-u)\|_{L^\infty([0,T_q], H^{-2})}=0.
	\end{equation}
	\textbf{Step 3: Prove the inequality \eqref{1.8}}.\\
	Let us now show that $u$ verifies inequality \eqref{1.8}. Let $t>0$ and let $n\in \N$, $R> 0$ and $q\in \N$ such that $t\in[0,T_q]$, we have
	\begin{align*}
		\|J_R(\theta_q u)(t)\|_{L^2}&\leq\|J_R			(\theta_q(u_{{\varphi(n)}}-u))(t)\|_{L^2}+\|J_R(\theta_q u_{{\varphi(n)}})(t)\|_{L^2}\\
		&\leq (2\pi)^{-\frac{3}{2}}(1+R^2)\|\theta_q(u_{{\varphi(n)}}-u)(t)\|_{L^2}+\|\theta_q u_{{\varphi(n)}}(t)\|_{L^2}\\
		&\leq(2\pi)^{-\frac{3}{2}}(1+R^2) \|\theta_q(u_{{\varphi(n)}}-u)\|_{L^\infty([0,T_q],H^{-2})} +\|u_{{\varphi(n)}}(t)\|_{L^2}.
	\end{align*}
	By passing to the $\liminf$ when $n$ tends to infinity, we get for all $R> 0$ and $q\in \N$
	\begin{align*}
		\|J_R(\theta_q u)(t)\|_{L^2}&\leq\liminf_{n\rightarrow \infty}\|u_{{\varphi(n)}}(t)\|_{L^2}.
	\end{align*}
	Due to the monotonic convergence theorem, when $R$ tends to infinity, we obtain for all $q\in \N$
	\begin{align*}
		\|\theta_q u(t)\|_{L^2}&\leq\liminf_{n\rightarrow \infty}\|u_{{\varphi(n)}}(t)\|_{L^2}.
	\end{align*}
	We applied the same theorem, when $q$ tends to infinity, we obtain
	\begin{equation}
		\label{3.5}	\|u(t)\|_{L^2}\leq \liminf_{n\rightarrow \infty}\|u_{{\varphi(n)}}(t)\|_{L^2}.
	\end{equation}
	In the other hand, we have the sequence $(u_{\varphi(n) })_{n}$ is bounded in the Hilbert space $L^2(\R^+,\dot{H}^1(\R^3)) $, then
	\begin{equation}
		u_{{\varphi(n)}}\rightarrow u {\;\rm weakly\; in\;} L^2( \R^+, \dot{H}^1(\R^3)).
	\end{equation}	
	Particularly, $u\in L^2(\R^+,\dot{H}^1(\R^3))$ and  for all $t\geq 0$
	\begin{equation}
		\label{3.7}	\int_{0}^{t}\|\nabla u\|_{L^2}^2\leq \liminf_{n\rightarrow \infty}\int_{0}^{t}\|\nabla u_{{\varphi(n)}}\|_{L^2}^2.
	\end{equation}		
	It now remains to prove that for all $t>0$
	\begin{equation}
		\int_0^t\int_{\R^3}f(|u|)|u|^2\leq \liminf_{n\rightarrow\infty}\int_0^t\int_{\R^3}f(|u_{{\varphi(n)}}|)|u_{{\varphi(n)}}|^2 .
	\end{equation}
	Let $n\in \N$ and $q\in \N$, then 
	\begin{align*}
		\|\theta_q(u_{\varphi(n)}-u)\|^2_{L^2([0,T_q]\times\R^3)}&=\int_{0}^{T_q}\|\theta_q(u_{\varphi(n)}(t)-u(t))\|_{L^2}^2dt\\
		&\leq\int_{0}^{T_q}\|\theta_q(u_{\varphi(n)}(t)-u(t))\|_{H^{-1}}\|\theta_q(u_{\varphi(n)}(t)-u(t))\|_{H^{1}}dt\\
		&\leq\left(\int_{0}^{T_q}\|\theta_q(u_{\varphi(n)}(t)-u(t))\|^2_{H^{-1}}dt\right)^{\frac{1}{2}}\times \left(\int_{0}^{T_q}\|\theta_q(u_{\varphi(n)}(t)-u(t))\|^2_{H^{1}}dt\right)^{\frac{1}{2}}\\
		&\leq\sqrt{T_q}\|\theta_q(u_{\varphi(n)}-u)\|_{L^\infty([0,T_q],H^{-1})}\times \left(\int_{0}^{T_q}\|\theta_q(u_{\varphi(n)}-u)\|^2_{H^{1}}dt\right)^{\frac{1}{2}}.
	\end{align*}
	In the other hand
	\begin{align*}
		\|\theta_q(u_{\varphi(n)}(t)-u(t))\|^2_{H^{1}}&=\|\theta_q(u_{\varphi(n)}(t)-u(t))\|^2_{L^2}+\|\nabla(\theta_q(u_{\varphi(n)}(t)-u(t)))\|^2_{L^2}\\
		&\leq\|u_{\varphi(n)}(t)-u(t)\|^2_{L^2}+\|(\nabla\theta_q)(u_{\varphi(n)}(t)-u(t))\|^2_{L^2}+\|\theta_q\nabla(u_{\varphi(n)}(t)-u(t))\|^2_{L^2}\\
		&\leq 2\|u^0\|_{L^2}^2+\|\nabla\theta_q\|^2_{L^\infty}\|u_{\varphi(n)}(t)-u(t)\|^2_{L^2}+\|\nabla u_{\varphi(n)}(t)\|^2_{L^2}+\|\nabla u(t)\|^2_{L^2}\\
		&\leq 2(1+\|\nabla\theta_q\|^2_{L^\infty})\|u^0\|_{L^2}^2+\|\nabla u_{\varphi(n)}(t)\|_{L^2}+\|\nabla u(t)\|^2_{L^2},
	\end{align*}
	which implies
	\begin{align*}
		\int_{0}^{T_q}\|\theta_q(u_{\varphi(n)}(t)-u(t))\|^2_{H^{1}}dt&\leq 2(1+\|\nabla\theta_q\|^2_{L^\infty})\int_{0}^{T_q}\|u^0\|_{L^2}^2dt+\int_{0}^{T_q}\|\nabla u_{\varphi(n)}(t)\|_{L^2}dt+\int_{0}^{T_q}\|\nabla u(t)\|^2_{L^2}dt\\
		&\leq (2T_q(1+\|\nabla\theta_q\|^2_{L^\infty})+2)\|u^0\|_{L^2}^2.
	\end{align*}
	By passing to the $\limsup$ when $n$ tends to infinity, we get
	\begin{align*}
		\limsup_{n\rightarrow \infty}\|\theta_q(u_{\varphi(n)}-u)\|^2_{L^2([0,T_q]\times\R^3)}&\leq\sqrt{T_q(2T_q(1+\|\nabla\theta_q\|^2_{L^\infty})+2)}\|u^0\|_{L^2}\limsup_{n\rightarrow \infty}\|\theta_q(u_{\varphi(n)}-u)\|_{L^\infty([0,T_q],H^{-1})}=0.
	\end{align*}
	Therefore, for all $q\in \N$, we have
	\begin{equation}
		\label{3.9}	\theta_qu_{\varphi(n)}\rightarrow \theta_q u{\;\;\;\rm strongly\; in\;} L^2( [0,T_q]\times\R^3). 
	\end{equation}
	For $q=1$, their exist a increasing function $\psi_1:\N \rightarrow \N $ such that
	\begin{equation}
		\theta_1u_{\varphi\circ \psi_1(n)}(t,x)\underset{n\rightarrow \infty}{\longrightarrow} \theta_1u(t,x) {\;\;\;\rm a.e.\;\;} (t,x)\in [0,T_1]\times \R^3,
	\end{equation}
	which implies that
	\begin{equation}
		u_{\varphi\circ \psi_1(n)}(t,x)\underset{n\rightarrow \infty}{\longrightarrow} u(t,x) {\;\;\;\rm a.e.\;\;} (t,x)\in [0,T_1]\times B(0,2).
	\end{equation}
	By recurrence, we construct a family of strictly increasing functions $(\psi_q)_{q\in \N}$ from $\N$ into $\N$ such that for all $q\in \N$
	\begin{equation}
		u_{\varphi\circ \psi_1\circ\cdots\circ\psi_q(n)}(t,x)\underset{n\rightarrow \infty}{\longrightarrow} u(t,x) {\;\;\;\rm a.e.\;\;} (t,x)\in [0,T_q]\times B(0,q+1).
	\end{equation}
	According to Cantor's diagonal process, there exists strictly increasing function $\psi:\N\rightarrow \N$ such that
	\begin{equation}
		u_{\varphi\circ \psi(n)}(t,x)\underset{n\rightarrow \infty}{\longrightarrow} u(t,x) {\;\;\;\rm a.e.\;\;} (t,x)\in \R^+\times \R^3.
	\end{equation}
	Since $f$ is a continuous function, we have 
	\begin{equation}
		f(|u_{\varphi\circ \psi(n)}(t,x)|)|u_{\varphi\circ \psi(n)}(t,x)|^2\underset{n\rightarrow \infty}{\longrightarrow} f(|u(t,x)|)|u(t,x)|^2 {\;\;\;\rm a.e.\;\;} (t,x)\in \R^+\times \R^3.
	\end{equation}
	Using Fatou's Lemma, we obtain for all $t>0$
	\begin{align}
		\label{3.15}	\int_0^t\int_{\R^3}f(|u(z,x)|)|u(z,x)|^2dxdz &\leq  \liminf_{n\rightarrow\infty} \int_0^t\int_{\R^3}f(|u_{{\varphi\circ\psi(n)}}(z,x)|)|u_{{\varphi\circ\psi(n)}}(z,x)|^2dxdz.
	\end{align}
	Let's note that  $\tilde{\varphi}=\varphi\circ \psi$, therefore
	\begin{align*}
		\eqref{3.5}\Rightarrow & \|u(t)\|_{L^2}\leq \liminf_{n\rightarrow \infty}\|u_{{\tilde{\varphi}(n)}}(t)\|_{L^2},\quad \forall t\geq 0,\\
		\eqref{3.7}\Rightarrow &\int_{0}^{t}\|\nabla u\|_{L^2}^2\leq \liminf_{n\rightarrow \infty}\int_{0}^{t}\|\nabla u_{{\tilde{\varphi}(n)}}\|_{L^2}^2,\quad \forall t\geq 0,\\
		\eqref{3.15}\Rightarrow &\int_0^t\int_{\R^3}f(|u|)|u|^2 \leq  \liminf_{n\rightarrow\infty} \int_0^t\int_{\R^3}f(|u_{{\tilde{\varphi}(n)}}|)|u_{{\tilde{\varphi}(n)}}|^2,\quad \forall t\geq 0.
	\end{align*}
	Combining the above inequalities with \eqref{Eq0}, we obtain: for all $t\geq 0$:
	\begin{equation}\label{3.16}
		\|u(t)\|^{2}_{L^{2}}+2\int_0^t\|\nabla u(s)\|^{2}_{L^{2}}ds+2\int_{0}^{t}\int_{\R^3}f(|u(s,x)|)|u(s,x)|^2dxds\leq \|u_{0}\|^{2}_{L^{2}}.
	\end{equation}
	%%%%%%%%%%%%%%%%%%%%%%%%%%%%%%%%%
	\textbf{Step 4: Show $u$ as a weak solution}.\\
	To show that $u=(u^1,u^2,u^3)$ is a weak solution of the \ref{NSD} equation, it is necessary to show that
	\begin{equation}
		\forall i,j\in\{1,2,3\},\quad J_{\tilde{\varphi}(n)}\big(u^i_{{\tilde{\varphi}(n)}}u^j_{{\tilde{\varphi}(n)}}\big)\underset{n\to+\infty}{\longrightarrow }u^iu^j,\mbox{ in }D'(\R^+\times \R^3), 
	\end{equation}
	and 
	\begin{equation}
		J_{\tilde{\varphi}(n)}\big(f(|u_{{\tilde{\varphi}(n)}}|)u^i_{{\tilde{\varphi}(n)}}\big)\underset{n\to+\infty}{\longrightarrow }f(|u|)
		u^i,\mbox{ in }D'(\R^+\times \R^3), 
	\end{equation}
	where $u_{{\tilde{\varphi}(n)}}=(u^1_{{\tilde{\varphi}(n)}},u^2_{{\tilde{\varphi}(n)}},u^3_{{\tilde{\varphi}(n)}})$.\\
	
	\noindent Let $\phi$ in $\mathcal{D}(\R^+\times \R^3)$. There is an integer $q\in \N$ such that
	$$\supp(\phi)\subset [0, T_q) \times B(0, q),$$
	and we have
	\begin{equation}
		\int_{0}^{+\infty}\int_{\R^3}\left[J_{\tilde{\varphi}(n)}\big(u^i_{{\tilde{\varphi}(n)}}u^j_{{\tilde{\varphi}(n)}}\big)-u^iu^j\right]\phi=I_n^1+I_n^2+I_n^3,
	\end{equation}
	where
	\begin{align*}
		&I_n^1=	-\int_{0}^{+\infty}\int_{\R^3}\left(Id-J_{\tilde{\varphi}(n)}\right)\big[u^i_{{\tilde{\varphi}(n)}}u^j_{{\tilde{\varphi}(n)}}\big]\phi;\\
		&I_n^2=	\int_{0}^{+\infty}\int_{\R^3}\big(u^i_{{\tilde{\varphi}(n)}}-u^i\big)u^j_{{\tilde{\varphi}(n)}}\phi;\\
		&I_n^3=	\int_{0}^{+\infty}\int_{\R^3}\big(u^j_{{\tilde{\varphi}(n)}}-u^j\big)u^i\phi.
	\end{align*}
	$\ast$ Estimate of term $I_n^1$: We have 
	\begin{align*}
		\left|I_n^1\right|&=\left|\int_{0}^{+\infty}\int_{\R^3}\left(Id-J_{\tilde{\varphi}(n)}\right)\big[u^i_{{\tilde{\varphi}(n)}}u^j_{{\tilde{\varphi}(n)}}\big]\phi\right|\\
		&\leq \int_{0}^{+\infty}\|\left(Id-J_{\tilde{\varphi}(n)}\right)u^i_{{\tilde{\varphi}(n)}}u^j_{{\tilde{\varphi}(n)}}\|_{\dot{H}^{-\frac{3}{2}}}\|\phi\|_{\dot{H}^{\frac{3}{2}}}\\
		&\leq \left(\int_{0}^{+\infty}\|\left(Id-J_{\tilde{\varphi}(n)}\right)u^i_{{\tilde{\varphi}(n)}}u^j_{{\tilde{\varphi}(n)}}\|_{\dot{H}^{-\frac{3}{2}}}^2\right)^{\frac{1}{2}}\|\phi\|_{L^2(\R^+,\dot{H}^{\frac{3}{2}})}\\
		&\leq \frac{1}{\tilde{\varphi}(n)}\left(\int_{0}^{+\infty}\|u^i_{{\tilde{\varphi}(n)}}u^j_{{\tilde{\varphi}(n)}}\|_{\dot{H}^{1-\frac{3}{2}}}^2\right)^{\frac{1}{2}}\|\phi\|_{L^2(\R^+,\dot{H}^{\frac{3}{2}})}\\
		&\leq \frac{C}{\tilde{\varphi}(n)}\left(\int_{0}^{+\infty}\|u^i_{{\tilde{\varphi}(n)}}\|^2_{L^2}\|\nabla u^j_{{\tilde{\varphi}(n)}}\|_{L^2}^2\right)^{\frac{1}{2}}\|\phi\|_{L^2(\R^+,\dot{H}^{\frac{3}{2}})}\\
		&\leq \frac{C}{\tilde{\varphi}(n)}\left(\int_{0}^{+\infty}\|u_{{\tilde{\varphi}(n)}}\|^2_{L^2}\|\nabla u_{{\tilde{\varphi}(n)}}\|_{L^2}^2\right)^{\frac{1}{2}}\|\phi\|_{L^2(\R^+,\dot{H}^{\frac{3}{2}})}\\
		&\leq \frac{C}{\tilde{\varphi}(n)}\|u^0\|_{L^2}\left(\int_{0}^{+\infty}\|\nabla u_{{\tilde{\varphi}(n)}}\|_{L^2}^2\right)^{\frac{1}{2}}\|\phi\|_{L^2(\R^+,\dot{H}^{\frac{3}{2}})}\\
		&\leq \frac{C}{\tilde{\varphi}(n)}\|u^0\|_{L^2}^2\|\phi\|_{L^2(\R^+,\dot{H}^{\frac{3}{2}})}.
	\end{align*}
	Therefore $$\lim\limits_{n\rightarrow +\infty}I_n^1=0.$$
	$\ast$ Estimate of term $I_n^2$: We have
	\begin{align*}
		|I_n^2|&=\left|\int_{0}^{+\infty}\int_{\R^3}\big(u^i_{{\tilde{\varphi}(n)}}-u^i\big)u^j_{{\tilde{\varphi}(n)}}\phi\right|\\
		&=\left|\int_{0}^{T_q}\int_{\R^3}\theta_q\big(u^i_{{\tilde{\varphi}(n)}}-u^i\big)u^j_{{\tilde{\varphi}(n)}}\phi\right|\\
		&\leq \|\theta_q\big(u^i_{{\tilde{\varphi}(n)}}-u^i\big)\|_{L^2([0,T_q]\times \R^3)}\|u^j_{{\tilde{\varphi}(n)}}\|_{L^2([0,T_q]\times \R^3)}\|\phi\|_{L^\infty(\R^+\times \R^3)}\\
		&\leq \sqrt{T_q}\|u^0\|_{L^2}\|\theta_q\big(u_{{\tilde{\varphi}(n)}}-u\big)\|_{L^2([0,T_q]\times \R^3)}\|\phi\|_{L^\infty(\R^+\times \R^3)}
	\end{align*}
	Using \eqref{3.9}, we obtain $$\lim\limits_{n\rightarrow \infty}I_n^2=0.$$
	In a similar way, we show that $$\lim\limits_{n\rightarrow \infty}I_n^3=0.$$
	Finally, we get
	\begin{equation}\label{Lim1}
		\lim_{n\rightarrow \infty} J_{\tilde{\varphi}(n)}\big(u^i_{{\tilde{\varphi}(n)}}u^j_{{\tilde{\varphi}(n)}}\big)=u^iu^j {\quad\rm in\quad} \mathcal{D}'(\R^+\times \R^3).
	\end{equation}
	
	\noindent In the other hand, we have
	\begin{align*}
		\int_{0}^{+\infty}\int_{\R^3}\left[J_{\tilde{\varphi}(n)}\left(f\left(|u_{\tilde{\varphi}(n)}|\right)u^i_{\tilde{\varphi}(n)}\right)-f\left(|u|\right)u^i\right]\phi\  =L^1_n+L^2_n+L^3_n,
	\end{align*}
	where
	\begin{align*}
		&L^1_n:=-\int_{0}^{+\infty}\int_{\R^3}\left(Id-J_{\tilde{\varphi}(n)}\right)\left[f\left(|u_{\tilde{\varphi}(n)}|\right)u^i_{\tilde{\varphi}(n)}\right]\phi\  ;\\
		&L^2_n:= \int_{0}^{+\infty}\int_{\R^3}f\left(|u|\right)\left[u^i_{\tilde{\varphi}(n)}-u^i\right]\phi\  ;\\
		&L^3_n:= \int_{0}^{+\infty}\int_{\R^3}u^i_{\tilde{\varphi}(n)}\left[f\left(|u_{\tilde{\varphi}(n)}|\right)-f\left(|u|\right)\right]\phi\ .
	\end{align*}
	We will show that for all $\varepsilon>0$, there exists a rank $N\in \N$ such that for all $n\geq N$
	$$|L_n^1|\leq \frac{\varepsilon}{3},\ |L_n^2|\leq \frac{\varepsilon}{3}\mbox{ and }|L_n^3|\leq \frac{\varepsilon}{3}.$$ For this we fix $\varepsilon>0$ and take $R=R_\varepsilon>0$ such that 
	$$\frac{1}{R}\leq \frac{\varepsilon}{24(\|u^0\|_{L^2}^2+1)\|\phi\|_{L^\infty(\R^+\times\R^3)}}.$$
	\begin{itemize}
		\item[$\ast$] \textbf{Estimate of the term $ L^1_n $:} We have 
		\begin{align*}
			\left|L^1_n\right|&=\left|\int_{0}^{T_q}\left(\left(Id-J_{\tilde{\varphi}(n)}\right)\left[f\left(|u_{\tilde{\varphi}(n)}|\right)u^i_{\tilde{\varphi}(n)}\right],\phi\right)_{L^2}\right|\\
			&\leq \frac{1}{\tilde{\varphi}(n)}\|f\left(|u_{\tilde{\varphi}(n)}|\right)u^i_{\tilde{\varphi}(n)}\|_{L^{1}([0,T_q],H^{-3})}\|\phi\|_{L^\infty(\R^+,H^4)}.
		\end{align*}
		Now, 
		\begin{align*}
			\|f\left(|u_{\tilde{\varphi}(n)}|\right)u^i_{\tilde{\varphi}(n)}\|_{L^{1}([0,T_q],H^{-3})}&\leq C\|f\left(|u_{\tilde{\varphi}(n)}|\right)u_{\tilde{\varphi}(n)}\|_{L^{1}([0,T_q]\times\R^3)}\\
			&\leq C \int_{\mathcal{A}_{n}}f\left(|u_{\tilde{\varphi}(n)}|\right)\left|u_{\tilde{\varphi}(n)}\right|dx+C\int_{\mathcal{A}_{n
				}^c}f\left(|u_{\tilde{\varphi}(n)}|\right)\left|u_{\tilde{\varphi}(n)}\right|dx,
		\end{align*}
		where 
		$$ \mathcal{A}_{n}:=\{(t,x)\in [0,T_q]\times\R^3;\ \left|u_{\tilde{\varphi}(n)}(t,x)\right|<1 \}.$$
		Using Lemma \ref{Lemma2.4}, for all $(t,x)\in \mathcal{A}_{n}$, we have
		$$f\left(|u_{\tilde{\varphi}(n)}(t,x)|\right)\left|u_{\tilde{\varphi}(n)}(t,x)\right|\leq f'(1)\left|u_{\tilde{\varphi}(n)}(t,x)\right|^2,$$
		which implies 
		$$\int_{\mathcal{A}_{n}}f\left(|u_{\tilde{\varphi}(n)}|\right)\left|u_{\tilde{\varphi}(n)}\right|dx\leq f'(1)\int_{0}^{T_q}\|u_{\tilde{\varphi}(n)}\|_{L^2}^2\leq T_qf'(1)\|u^0\|_{L^2}^2.$$
		Moreover, for all $x\in  \mathcal{A}_{n}^c$, we have
		$$f\left(|u_{\tilde{\varphi}(n)}(t,x)|\right)\left|u_{\tilde{\varphi}(n)}(t,x)\right|\leq f\left(|u_{\tilde{\varphi}(n)}(t,x)|\right)\left|u_{\tilde{\varphi}(n)}(t,x)\right|^2,$$
		which implies 
		$$\int_{\mathcal{A}_{n}^c}f\left(|u_{\tilde{\varphi}(n)}|\right)\left|u_{\tilde{\varphi}(n)}\right|dx\leq\|f\left(|u_{\tilde{\varphi}(n)}|\right)\left|u_{\tilde{\varphi}(n)}\right|^2\|_{L^1(\R^+\times \R^3)}\leq \|u^0\|_{L^2}^2.$$
		Thus,
		\begin{align*}
			\|f\left(|u_{\tilde{\varphi}(n)}|\right)u^i_{\tilde{\varphi}(n)}\|_{L^{1}([0,T_q],H^{-3})}&\leq C \left(T_q f'(1)+1\right)\|u^0\|_{L^2}^2.
		\end{align*}
		As a result, 
		\begin{align*}
			|L^1_n|&\leq \frac{C}{\tilde{\varphi}(n)}\left(T_qf'(1)+1\right)\|u^0\|_{L^2}^2\|\phi\|_{L^\infty(\R^+,H^4)}.
		\end{align*}
		This implies that there exists an integer $ n_1 \in \N $ such that for all $ n \geq n_1$,
		$$	|L^1_n|\leq\frac{\varepsilon}{3}.$$
		\item[$\ast$] \textbf{Estimate of the term $ L^2_n $:} We have  
		\begin{align*}
			\left|L^2_n\right|&\leq  \int_{\mathcal{A}}f\left(|u|\right)\left|u_{\tilde{\varphi}(n)}-u\right|\left|\phi\right|+\int_{\mathcal{A}^c}f\left(|u|\right)\left|u_{\tilde{\varphi}(n)}-u\right|\left|\phi\right|
		\end{align*}
		where	$\mathcal{A}=\{(t,x)\in \R^+\times\R^3;\ |u(t,x)|\leq R\}.$ \\
		Using the fact that $ f $ is an increasing function, we obtain
		\begin{align*}
			\int_{\mathcal{A}}f\left(|u|\right)\left|u_{\tilde{\varphi}(n)}-u\right|\left|\phi\right|&\leq f\left(R\right)\int_{0}^{T_q}\int_{\R^3}\left|\theta_q(u_{\tilde{\varphi}(n)}-u)\right|\left|\phi\right|\\
			&\leq f(R)\|\theta_q(u_{\tilde{\varphi}(n)}-u)\|_{L^2([0,T_q]\times \R^3)}\|\phi\|_{L^2(\R^+\times\R^3)}.
		\end{align*}
		Using \eqref{3.9}, we deduce that there exists $ n_2 \in \N $ such that for all $n\geq n_2$,
		\begin{equation}
			\int_{\mathcal{A}}f\left(|u|\right)\left|u_{\tilde{\varphi}(n)}-u\right|\left|\phi\right|\leq\frac{\varepsilon}{12}.
		\end{equation}
		On the other hand, define 
		\begin{equation}\label{BN}
			\mathcal{B}_n:=\left\{(t,x)\in \R^+\times \R^3;\ |u(t,x)|\leq |u_{\tilde{\varphi}(n)}(t,x)|\right\}.
		\end{equation}
		Then,
		\begin{align*}
			\int_{\mathcal{A}^c}f\left(|u|\right)\left|u_{\tilde{\varphi}(n)}-u\right|\left|\phi\right|&\leq \int_{\mathcal{A}^c}f\left(|u|\right)\left|u_{\tilde{\varphi}(n)}\right|\left|\phi\right|+\int_{\mathcal{A}^c}f\left(|u|\right)\left|u\right|\left|\phi\right|\\
			&\leq\int_{\mathcal{A}^c\cap\mathcal{B}_n}f\left(|u|\right)\left|u_{\tilde{\varphi}(n)}\right|\left|\phi\right|+2\int_{\mathcal{A}^c}f\left(|u|\right)\left|u\right|\left|\phi\right|.
		\end{align*}
		Use the fact that $f$ es an increasing function, we get
		\begin{align*}
			\int_{\mathcal{A}^c\cap\mathcal{B}_n}f\left(|u|\right)\left|u_{\tilde{\varphi}(n)}\right|\left|\phi\right|&\leq \int_{\mathcal{A}^c\cap\mathcal{B}_n}f\left(|u_{\tilde{\varphi}(n)}|\right)\left|u_{\tilde{\varphi}(n)}\right|\left|\phi\right|\\
			&\leq \frac{1}{R} \int_{\mathcal{A}^c\cap\mathcal{B}_n}f\left(|u_{\tilde{\varphi}(n)}|\right)\left|u_{\tilde{\varphi}(n)}\right|^2\left|\phi\right|\\
			&\leq \frac{1}{R}\|\phi\|_{L^\infty(\R^+\times\R^3)}\|f\left(|u_{\tilde{\varphi}(n)}|\right){\left|u_{\tilde{\varphi}(n)}\right|^2}\|_{L^1(\R^+\times\R^3)}\\
			&\leq \frac{1}{R}\|\phi\|_{L^\infty(\R^+\times\R^3)}\|u^0\|_{L^2}^2\\
			&\leq \frac{\varepsilon}{12}.
		\end{align*}
		Similarly,
		\begin{align*}
			\int_{\mathcal{A}^c}f\left(|u|\right)\left|u\right|\left|\phi\right|\leq \frac{1}{R} \|\phi\|_{L^\infty(\R^+\times\R^3)}\|u^0\|_{L^2}^2\leq \frac{\varepsilon}{12}.
		\end{align*}
		Combining the above results, we conclude that there exists $ n_2 \in \N $ such that for all $ n \geq n_2 $,
		$$	|L^2_n|\leq\frac{\varepsilon}{3}.$$
		\item[$\ast$] \textbf{Estimate of the term $ L^3_n $:} We have
		\begin{align*}
			\left|L^3_n\right|&\leq  \int_{\mathcal{A}_{n,R}}\left|u_{\tilde{\varphi}(n)}\right|\left|f\left(|u_{\tilde{\varphi}(n)}|\right)-f\left(|u|\right)\right|\left|\phi\right|+ \int_{\mathcal{A}_{n,R}^c}\left|u_{\tilde{\varphi}(n)}\right|\left|f\left(|u_{\tilde{\varphi}(n)}|\right)-f\left(|u|\right)\right|\left|\phi\right|:=L^{3,1}_n+L^{3,2}_n,
		\end{align*}
		where
		$$\mathcal{A}_{n,R}:=\left\{(t,x)\in \R^+\times\R^3;\ |u_{\tilde{\varphi}(n)}(t,x)|\leq R\right\}.$$
		Now,
		\begin{align*}
			L^{3,1}_n&=\int_{\mathcal{A}_{n,R}\cap\mathcal{B}_{n}}\left|u_{\tilde{\varphi}(n)}\right|\left|f\left(|u_{\tilde{\varphi}(n)}|\right)-f\left(|u|\right)\right|\left|\phi\right|+\int_{\mathcal{A}_{n,R}\cap\mathcal{B}_{n}^c}\left|u_{\tilde{\varphi}(n)}\right|\left|f\left(|u_{\tilde{\varphi}(n)}|\right)-f\left(|u|\right)\right|\left|\phi\right|,
		\end{align*}
		where $\mathcal{B}_n$ is defined in \eqref{BN}.\\
		
		Let $(t,x)\in \mathcal{A}_{n,R}\cap\mathcal{B}_{n}$. According to Lemma \ref{Lemma2.4}, we have
		$$\left|u_{\tilde{\varphi}(n)}(t,x)\right|\left|f\left(|u_{\tilde{\varphi}(n)}(t,x)|\right)-f\left(|u(t,x)|\right)\right|\leq R f'(R)\left|u_{\tilde{\varphi}(n)}(t,x)-u(t,x)\right|.$$
		Thus,
		\begin{align*}
			\int_{\mathcal{B}_{n}\cap\mathcal{A}_{n,R}}\left|u_{\tilde{\varphi}(n)}\right|\left|f\left(|u_{\tilde{\varphi}(n)}|\right)-f\left(|u|\right)\right|\left|\phi\right|&\leq Rf'(R)\int_{0}^{T_q}\int_{\R^3}\left|\theta_q(u_{\tilde{\varphi}(n)}-u)\right|\left|\phi\right|\\
			&\leq Rf'(R)\|\phi\|_{L^{2}(\R^+\times\R)}\|\theta_q(u_{\tilde{\varphi}(n)}-u)\|_{L^2([0,T_q]\times \R^3)}.
		\end{align*}
		This implies that there exists an integer $ n_3 \in \N$ such that, for all $n \geq n_3 $,
		$$	\int_{\mathcal{A}_{n,R}\cap \mathcal{B}_n}\left|u_{\tilde{\varphi}(n)}\right|\left|f\left(|u_{\tilde{\varphi}(n)}|\right)-f\left(|u|\right)\right|\left|\phi\right|\leq\frac{\varepsilon}{12}.$$
		On the other hand, for $(t,x)\in \mathcal{A}_{n,R}\cap \mathcal{B}_n^c$, we have 
		\begin{align*}\left|u_{\tilde{\varphi}(n)}(t,x)\right|\left|f\left(|u_{\tilde{\varphi}(n)(t,x)}|\right)-f\left(|u(t,x)|\right)\right|&\leq\left|u(t,x)\right|\left|f\left(|u_{\tilde{\varphi}(n)(t,x)}|\right)-f\left(|u(t,x)|\right)\right|\\
			&\leq \begin{cases}
				Rf'(R) \left|u_{\tilde{\varphi}(n)}(t,x)-u(t,x)\right|&\mbox{ if }(t,x)\in \mathcal{A}_R\\ \\
				\frac{1}{R}f(|u(t,x)|)|u(t,x)|^2&\mbox{ if }(t,x)\in \mathcal{A}_R^c
			\end{cases}\\
			&\leq 	Rf'(R) \left|u_{\tilde{\varphi}(n)(t,x)}-u(t,x)\right|+\frac{1}{R}f(|u(t,x)|)|u(t,x)|^2,
		\end{align*}
		where
		$$\mathcal{A}_{R}:=\left\{(t,x)\in \R^+\times\R^3;\ |u(t,x)|\leq R\right\}.$$
		Therefore,
		\begin{align*}
			&	\hskip-1cm\int_{\mathcal{A}_{n,R}\cap\mathcal{B}_{n}^c}\left|u_{\tilde{\varphi}(n)}\right|\left|f\left(|u_{\tilde{\varphi}(n)}|\right)-f\left(|u|\right)\right|\left|\phi\right|\\
			&\leq Rf'(R)\int_{0}^{T_q}\int_{\R^3}\left|\theta_q(u_{\tilde{\varphi}(n)}-u)\right|\left|\phi\right|+ \frac{1}{R}\int_{0}^{T_q}\int_{\R^3}f(|u|)\left|u\right|^2\left|\phi\right|\\
			&\leq Rf'(R)\|\theta_q(u_{\tilde{\varphi}(n)}-u)\|_{L^{2}([0,T_q]\times\R^3)} \|\phi\|_{L^2(\R^+\times \R)}+\frac{1}{R}\|u^0\|^2_{L^2}\|\phi\|_{L^\infty(\R^+\times \R)}\\
			&\leq Rf'(R)\|\theta_q(u_{\tilde{\varphi}(n)}-u)\|_{L^{2}([0,T_q]\times\R^3)} \|\phi\|_{L^2(\R^+\times \R)}+\frac{\varepsilon}{12}.
		\end{align*}
		This implies that there exists an integer $n_4\in \N$ such that, for all $n\geq n_4$,
		$$	\int_{\mathcal{A}_{n,R}\cap \mathcal{B}_n^c}\left|u_{\tilde{\varphi}(n)}\right|\left|f\left(|u_{\tilde{\varphi}(n)}|\right)-f\left(|u|\right)\right|\left|\phi\right|\leq\frac{\varepsilon}{6}.$$
		Let $n_5=n_3+n_4$. For all $n\geq n_5$,
		$$|L_n^{3,1}|\leq \dfrac{\varepsilon}{4}.$$
		Finally, for \( L_n^{3,2} \), using the fact that \( f \) is an increasing function, we have, for $(t,x)\in \mathcal{A}_{n,R}$,
		\begin{align*}
			\left|u_{\tilde{\varphi}(n)}(t,x)\right|\left|f\left(|u_{\tilde{\varphi}(n)}(t,x)|\right)-f\left(|u(t,x)|\right)\right|&\leq \begin{cases}
				f\left(|u_{\tilde{\varphi}(n)}(t,x)|\right)	\left|u_{\tilde{\varphi}(n)}(t,x)\right|&\mbox{if }(t,x)\in \mathcal{B}_n\\ \\
				f\left(|u(t,x)|\right)	\left|u(t,x)\right|&\mbox{if }(t,x)\in \mathcal{B}_n^c
			\end{cases}\\
			&\leq \begin{cases}
				\frac{1}{R}f\left(|u_{\tilde{\varphi}(n)}(t,x)|\right)	\left|u_{\tilde{\varphi}(n)}(t,x)\right|^2&\mbox{if }(t,x)\in \mathcal{B}_n\\ \\
				\frac{1}{R}	f\left(|u(t,x)|\right)	\left|u(t,x)\right|^2&\mbox{if }(t,x)\in \mathcal{B}_n^c.
			\end{cases}
		\end{align*}
		Thus,
		\begin{equation}
			\left|u_{\tilde{\varphi}(n)}\right|\left|f\left(|u_{\tilde{\varphi}(n)}|\right)-f\left(|u|\right)\right|\leq \frac{1}{R}\left(f\left(|u_{\tilde{\varphi}(n)}|\right)	\left|u_{\tilde{\varphi}(n)}\right|^2+f\left(|u|\right)	\left|u\right|^2\right).
		\end{equation}
		Therefore
		\begin{align*}
			|L_n^{3,2}|&\leq\frac{1}{R} \left[\int_{\R^+\times \R^3}f\left(|u_{\tilde{\varphi}(n)}|\right)	\left|u_{\tilde{\varphi}(n)}\right|^2|\phi|+\int_{\R^+\times \R^3}f\left(|u|\right)	\left|u\right|^2|\phi|\right]\\
			&\leq \frac{2}{R}\|u^0\|_{L^2}^2\|\phi\|_{L^2(\R^+\times\R^3)}\\
			&\leq \frac{\varepsilon}{12}.
		\end{align*}
		Finally,  there exists $n_5\in \N$ such that, for all $n\geq n_5$,
		$$\left|L^3_n\right|\leq \frac{\varepsilon}{3}.$$
	\end{itemize}
	As a result, we assume $N=n_1+n_2+n_3+n_4$, then for any $n\geq N$ we have 
	\begin{align*}
		\left|\int_{0}^{+\infty}\int_{\R^3}\left[J_{\tilde{\varphi}(n)}\left(f\left(|u_{\tilde{\varphi}(n)}|\right)u_{\tilde{\varphi}(n)}\right)-f\left(|u|\right)u\right]\phi\  \right|\leq \left|L^1_n\right|+\left|L^2_n\right|+\left|L^3_n\right|<\varepsilon.
	\end{align*}
	\subsection{Uniqueness of Solutions} Let $(u,\pi_u)$ and $ (v,\pi_v)$ be two weak solutions of \ref{NSD}. Then,
	\begin{align}
		\begin{cases}
			\;\partial_t u-\nu \Delta u+u\cdot\nabla u +f(|u|)u&=-\nabla \pi_u \\
			\;\partial_t v-\nu \Delta v+v\cdot\nabla v +f(|v|)v&=-\nabla \pi_v
		\end{cases}
	\end{align}
	Define $\omega=u-v$. Then $\omega$ satisfies the following equation:
	\begin{equation}\label{3.22}
		\partial_t \omega -\nu\Delta \omega+\omega\cdot\nabla u+v\cdot\nabla \omega +(f(|u|)u-f(|v|)v)=-\nabla(\pi_u-\pi_v)
	\end{equation}
	Taking the $ L^2 $-inner product of \eqref{3.22} with $ \omega $ and using the fact that $ \dive(u) = \dive(v) = 0 $, we obtain:
	\begin{align*}
		\frac{1}{2}\frac{d}{dt}\|w\|^{2}_{L^2}+\nu\|\nabla w\|^{2}_{L^2} +\left(f(|u|)u-f(|v|)v,w\right)_{L^2}&=\left(\dive(\omega \otimes u),\omega\right)_{L^2}\\
		&\leq \|\omega u\|_{L^2}\|\nabla \omega\|_{L^2}\\
		&\leq \frac{1}{4\nu}\|\omega u\|_{L^2}^2+\nu\|\nabla \omega\|_{L^2}^2
	\end{align*}
	This implies:
	\begin{align*}
		\frac{1}{2}\frac{d}{dt}\|w\|^{2}_{L^2} +\left(f(|u|)u-f(|v|)v,w\right)_{L^2}\leq \frac{1}{4\nu}\|\omega u\|_{L^2}^2.
	\end{align*}
	Using Lemma \ref{Lemma2.5}, for all $(t,x)\in \R^+\times \R^3$, we have:
	\begin{align*}
		\left(f(|u|)u-f(|v|)v,w\right)_{\R^3}&=\left(f(|u|)u-f(|v|)v,u-v\right)_{\R^3}\\
		&\geq \frac{c}{4}\left(|u|^{p}+|v|^{p}\right)|u-v|^2\\
		&\geq \frac{c}{4}|u|^{p}|\omega|^2.
	\end{align*}
	Thus:
	\begin{align*}
		\left(f(|u|)u-f(|v|)v,w\right)_{L^2}&=\int_{\R^3}\left(f(|u|)u-f(|v|)v,w\right)_{\R^3}\\
		&\geq \frac{c}{2}\int_{\R^3}|u|^{p} \left|\omega\right|^2
	\end{align*}
	On the other hand, using the fact that $ p > 1 $, we have:
	\begin{align*}
		\dfrac{1}{4\nu}	\int_{\R^3} |\omega|^2|u|^2&\leq c_{\nu,p}\int_{\R^3}|\omega|^2+\frac{c}{2}\int_{\R^3}|u|^{p} \left|\omega\right|^2.
	\end{align*}
	This implies:
	\begin{align*}
		\frac{1}{2}\frac{d}{dt}\|w\|^{2}_{L^2} +\frac{c}{2}\int_{\R^3}|u|^{p} \left|\omega\right|^2\leq  c_{\nu,p}\int_{\R^3}|\omega|^2+\frac{c}{2}\int_{\R^3}|u|^{p} \left|\omega\right|^2.
	\end{align*}
	Simplifying, we get:
	\begin{align*}
		\frac{d}{dt}\|w(t)\|^{2}_{L^2}
		&\leq c_{\nu,p}\|\omega(t)\|^2
	\end{align*} 
	Since $\omega(0)=0$, applying Grönwall's Lemma yields $\omega(t,x)=0$, for all $(t,x)\in \R^+\times \R^3$.\\
	
	Thus $u=v$, proving the uniqueness of the solution.
	\subsection{Continuity} The proof of the continuity of the weak solution $u$ is inspired by the proof given by Blel and Benameur in \cite{BB1}.
	We apply the upper limit in inequality \eqref{3.16} as $ t \to 0 $, obtaining
	\begin{equation}
		\limsup_{t\rightarrow 0}\|u(t)\|_{L^2}\leq \|u^0\|_{L^2}.
	\end{equation}
	Using the fact that 
	\begin{equation}
		\lim\limits_{t\rightarrow 0}\left(u(t)-u^0,v\right)_{L^2}=0, \quad \forall v \in L^2(\R^3),
	\end{equation}
	and according to Proposition \ref{prop1}, we have
	\begin{equation}
		\lim\limits_{t\rightarrow 0}\|u(t)-u^0\|_{L^2}=0,
	\end{equation}
	which implies that $ u $ is continuous at $ 0 $ in $ L^2(\R^3) $.\\
	Let $ t > 0 $ and $ \varepsilon \in (0, t) $, and define the functions $ v_\varepsilon = u(t \pm \varepsilon) $ and $ \pi_\varepsilon = \pi(t \pm \varepsilon) $. Then,
	\begin{align}
		\begin{cases}
			\;\partial_t u-\nu \Delta u+u\cdot\nabla u +f(|u|)u=-\nabla \pi,\\
			\;\partial_t v_\varepsilon-\nu \Delta v_\varepsilon+v_\varepsilon\cdot\nabla v_\varepsilon +f(|v_\varepsilon|)v_\varepsilon=-\nabla \pi_\varepsilon.
		\end{cases}
	\end{align}
	As a result, if $ \omega_\varepsilon = v_\varepsilon - u $, then
	\begin{equation}
		\partial_t \omega_{\varepsilon} -\nu\Delta \omega_\varepsilon+v_\varepsilon\cdot\nabla \omega_\varepsilon+\omega_\varepsilon\cdot\nabla u +f(|v_\varepsilon|)v_\varepsilon-f(|u|)u=-\nabla( \pi_\varepsilon-\pi).
	\end{equation}
	We apply the scalar product with $ \omega_\varepsilon $ in $ L^2(\R^3) $ and use the fact that $ \mathrm{div}(\omega_\varepsilon) = 0 $, obtaining
	\begin{align*}
		\frac{1}{2}\frac{d}{dt}\|\omega_\varepsilon(t)\|_{L^2}^2+\nu\|\nabla\omega_\varepsilon(t)\|_{L^2}^2+\left(f(|v_\varepsilon|)v_\varepsilon-f(|u|)u,\omega_\varepsilon\right)_{L^2}=-\left(\omega_\varepsilon\cdot\nabla u,\omega_\varepsilon\right)_{L^2}.
	\end{align*}
	Using the same technique as in the uniqueness proof, we obtain
	\begin{align*}
		\left(f(|v_\varepsilon|)v_\varepsilon-f(|u|)u,\omega_\varepsilon\right)_{L^2}\geq \frac{c}{2}\int_{\R^3}|u|^{p}|\omega_\varepsilon|^2,
	\end{align*}
	and 
	\begin{align*}
		\left|\left(\omega_\varepsilon\cdot\nabla u,\omega_\varepsilon\right)_{L^2}\right|\leq \nu \|\nabla\omega_\varepsilon\|_{L^2}^2+\frac{c}{2}\int_{\R^3}|u|^{p}|\omega_\varepsilon|^2+c_{\nu,p}\|\omega_\varepsilon\|_{L^2}^2.
	\end{align*}
	This implies that 
	\begin{align*}
		\frac{d}{dt}\|\omega_\varepsilon(t)\|_{L^2}^2\leq 2c_{\nu,p}\|\omega_\varepsilon(t)\|_{L^2}^2.
	\end{align*}
	Applying integration over $ (\varepsilon, t) $, we obtain
	\begin{align*}
		\|\omega_\varepsilon(t)\|_{L^2}^2\leq \|\omega_\varepsilon(\varepsilon)\|_{L^2}^2+\int_{\varepsilon}^{t} 2c_{\nu,p}\|\omega_\varepsilon(z)\|_{L^2}^2dz.
	\end{align*}
	Gronwall's Lemma then gives
	\begin{align*}
		\|\omega_\varepsilon(t)\|_{L^2}^2\leq \|\omega_\varepsilon(\varepsilon)\|_{L^2}^2 e^{2c_{\nu,p}(t-\varepsilon)},
	\end{align*}
	which implies
	\begin{align*}
		\|v_\varepsilon(t)-u(t)\|_{L^2}\leq \|v_\varepsilon(\varepsilon)-u(\varepsilon)\|_{L^2}e^{c_{\nu,p}t}.
	\end{align*}
	
	There are two cases to consider for $ v_\varepsilon = u(t+\varepsilon) $ and $ v_\varepsilon = u(t-\varepsilon) $:
	\begin{enumerate}
		\item[$\bullet$] If $ v_\varepsilon = u(t+\varepsilon) $, then 
		\begin{align*}
			\|u(t+ \varepsilon)-u(t)\|_{L^2}&\leq \|u(2\varepsilon)-u(\varepsilon)\|_{L^2}e^{c_{\nu,p}t}\\
			&\leq  \left(\|u(2\varepsilon)-u^0\|_{L^2}+\|u(\varepsilon)-u^0\|_{L^2}\right)e^{c_{\nu,p}t}.
		\end{align*}
		Taking $ \ds\limsup_{\varepsilon \to 0^+} $, we obtain
		\begin{align*}
			\limsup_{\varepsilon \to 0^+}\|u(t+ \varepsilon)-u(t)\|_{L^2}&\leq  \left(\limsup_{\varepsilon \to 0^+}\|u(2\varepsilon)-u^0\|_{L^2}+\limsup_{\varepsilon \to 0^+}\|u(\varepsilon)-u^0\|_{L^2}\right)e^{c_{\nu,p}t}=0.
		\end{align*}
		Thus,
		\begin{equation}
			\lim\limits_{\varepsilon \to 0^+}\|u(t+ \varepsilon)-u(t)\|_{L^2}=0,
		\end{equation}
		which implies that $ u $ is right-continuous at every point in $ \R^+ $.
		
		\item[$\bullet$] If $ v_\varepsilon = u(t-\varepsilon) $, then 
		\begin{align*}
			\|u(t- \varepsilon)-u(t)\|_{L^2}&\leq \|u(0)-u(\varepsilon)\|_{L^2}e^{c_{\nu,p}t}= \|u(\varepsilon)-u^0\|_{L^2}e^{c_{\nu,p}t}.
		\end{align*}
		Taking $\ds \limsup_{\varepsilon \to 0^+} $, we obtain
		\begin{align*}
			\limsup_{\varepsilon \to 0^+}\|u(t- \varepsilon)-u(t)\|_{L^2}&\leq  \limsup_{\varepsilon \to 0^+}\|u(\varepsilon)-u^0\|_{L^2}e^{c_{\nu,p}t}=0.
		\end{align*}
		Thus,
		\begin{equation}
			\lim\limits_{\varepsilon \to 0^+}\|u(t-\varepsilon)-u(t)\|_{L^2}=0,
		\end{equation}
		which implies that $ u $ is left-continuous at every point in $ \R^{+,\ast} $.	
	\end{enumerate}
	
	Therefore, $ u \in C_b(\R^+, L^2(\R^3)) $.\hfill$\blacksquare$
	\vskip0.5cm
	\textbf{Data availability} The data are available from the corresponding author upon reasonable request.
\section*{\bf Declarations}
	\textbf{Conflict of interest} The authors declare that there is no conflict of interest


\medskip
	
	
	\begin{thebibliography}{10}
		\bibitem{HB} BAHOURI, H. Fourier Analysis and Nonlinear Partial Differential Equations. Grundlehren der Mathematischen Wissenschaften, 2011, vol. 343.
		
		\bibitem{JBB} BENAMEUR, J. Global weak solution of 3D-NSE with exponential damping. Open Mathematics, 2022, vol. 20, no 1, p. 590-607.
		
		\bibitem{BB} BLEL, M. and BENAMEUR, J. Asymptotic analysis of Leray solution for the incompressible NSE with damping. Demonstratio Mathematica, 2024, vol. 57, no 1, p. 20240042.
		
		\bibitem{BB1} BLEL, M. and BENAMEUR, J. Long-time decay of leray solution of 3d-nse with exponential damping. Fractals, 2022, vol. 30, no 10, p. 2240236.
		
		\bibitem{HBAF} BREZIS, H. Analyse fonctionnelle Th eorie et application, 3e tirage ed. 1992.
		
		\bibitem{CJ} CAI, X and JIU, Q. Weak and strong solutions for the incompressible Navier–Stokes equations with damping. Journal of Mathematical Analysis and Applications, 2008, vol. 343, no 2, p. 799-809.
		
		\bibitem{JYC} CHEMIN, J-Y. DESJARDINS, Benoît, GALLAGHER, Isabelle, et al. Fluids with anisotropic viscosity. ESAIM: Mathematical Modelling and Numerical Analysis, 2000, vol. 34, no 2, p. 315-335.
		
		\bibitem{H} HOPF, E. Über die Anfangswertaufgabe für die hydrodynamischen Grundgleichungen. Erhard Schmidt zu seinem 75. Geburtstag gewidmet. Mathematische Nachrichten, 1950, vol. 4, no 1-6, p. 213-231.
		
		\bibitem{HL} HSIAO, L. Quasilinear hyperbolic systems and dissipative mechanisms. World Scientific, 1997.
		
		\bibitem{HP} HUANG, F. and PAN, R. Convergence rate for compressible Euler equations with damping and vacuum. Archive for Rational Mechanics and Analysis, 2003, vol. 166, p. 359-376.
		
		
		\bibitem{L} LERAY, J. Sur le mouvement d'un liquide visqueux emplissant l'espace. Acta mathematica, 1934, vol. 63, p. 193-248.
		
		
		
		
	\end{thebibliography}
\end{document}